# BUMP HUNTING WITH NON-GAUSSIAN KERNELS

By Peter Hall, Michael C. Minnotte and Chunming Zhang[1]

*Australian National University, Utah State University
and University of Wisconsin*

It is well known that the number of modes of a kernel density estimator is monotone nonincreasing in the bandwidth if the kernel is a Gaussian density. There is numerical evidence of nonmonotonicity in the case of some non-Gaussian kernels, but little additional information is available. The present paper provides theoretical and numerical descriptions of the extent to which the number of modes is a nonmonotone function of bandwidth in the case of general compactly supported densities. Our results address popular kernels used in practice, for example, the Epanechnikov, biweight and triweight kernels, and show that in such cases nonmonotonicity is present with strictly positive probability for all sample sizes $n \geq 3$. In the Epanechnikov and biweight cases the probability of nonmonotonicity equals 1 for all $n \geq 2$. Nevertheless, in spite of the prevalence of lack of monotonicity revealed by these results, it is shown that the notion of a critical bandwidth (the smallest bandwidth above which the number of modes is guaranteed to be monotone) is still well defined. Moreover, just as in the Gaussian case, the critical bandwidth is of the same size as the bandwidth that minimises mean squared error of the density estimator. These theoretical results, and new numerical evidence, show that the main effects of nonmonotonicity occur for relatively small bandwidths, and have negligible impact on many aspects of bump hunting.

**1. Introduction.** Compactly supported kernels, particularly the biweight, are predominantly used in practice when constructing a nonparametric density estimator. There are at least two reasons: ease of computation (calculation is simplified if a curve estimate at a given point uses only a relatively

Received January 2002; revised August 2003.
[1]Supported in part by NSF Grant DMS-03-53941.
*AMS 2000 subject classifications.* Primary 62G07; secondary 62G20.
*Key words and phrases.* Bandwidth choice, bootstrap, critical bandwidth, density estimation, kernel methods, modality, mode test, nonparametric curve estimation, unimodality.

small fraction of the data); and, more philosophically, a desire to ensure that a density estimator uses only local information. However, many "shape" properties of kernel density estimators are well understood only in the case of infinitely supported, Gaussian kernels.

Responding to this issue, in the present paper we quantify a range of properties of non-Gaussian kernels when used to identify bumps in nonparametric density estimation. Numerical results [e.g., Minnotte and Scott (1993)] have shown that in practical circumstances some commonly used, compactly supported kernels may not give rise to nonparametric density estimators whose modality is a monotone function of bandwidth. However, theoretical explanations of this property are not available, and neither is it clear whether the nonmonotonicity property will upset the sorts of bump-hunting applications to which kernel density estimators are often put. For example, can the biweight kernel be profitably used to implement Silverman's (1981) test for unimodality, or does the nonmonotonicity property interfere at too high a level for this to be feasible?

This paper provides answers to these questions. In Section 2 we give general theoretical results that address nonmonotonicity problems arising with compactly supported kernels. The results are illustrated theoretically in terms of commonly employed "multiweight" kernels, such as the uniweight (or Epanechnikov) density, and the biweight and triweight densities. Nevertheless our results are very general, and apply to a wide range of compactly supported kernels. Numerical illustrations of theoretical properties are given in Section 4.

To give an example, it follows from our results that when density estimators are calculated using the biweight kernel, the number of modes of a kernel density estimator is, with probability 1, a nonmonotone function of the bandwidth, $h$, whenever sample size, $n$, equals two or more. Interestingly, this result fails for the triweight kernel. In that case, for $n = 2$ and with probability 1, the number of modes is monotone in $h$. However the probability that it is nonmonotone is strictly positive whenever $n \geq 3$.

Results of this type add considerably to the information provided by more conventional analytical results, such as those of Schoenberg (1950). From those it may be deduced only that for compactly supported kernels, and sufficiently large sample sizes, there exist deterministic data constructions for which the nonmonotonicity property fails. By way of contrast, our results show that nonmonotonicity fails for the sorts of datasets that arise in practice, and for a wide range of sample sizes (generally, for $n \geq 3$).

These properties lead pointedly to the question of whether the critical bandwidth, in the case of compactly supported kernels, is of the same size as it would be for a Gaussian kernel. The critical bandwidth is defined as the "smallest" bandwidth, in some sense, such that a nonparametric density estimator is unimodal. When the Gaussian kernel is used, monotonicity of the number of modes as a function of bandwidth means that there is no ambiguity in the definition of "smallest."

The situation is much less clear for compactly supported kernels, however. Nevertheless, we shall show that provided the kernel is unimodal and concave at the mode, one may unambiguously define the "smallest" bandwidth $h_{\text{crit}}$ to be the infimum of values $h_1 > 0$ such that the number of modes of a density estimator equals 1 for all $h \geq h_1$. This version of $h_{\text{crit}}$ is well defined, and strictly positive, with probability 1. Moreover, the mode tree technology developed by Minnotte and Scott (1993) [see also Minnotte (1997)] enables this definition to be used in practice without difficulty.

One of the particularly attractive features of the critical bandwidth for a Gaussian kernel is that it is of size $n^{-1/5}$, this being the order that produces optimal mean squared error performance for density estimators in a standard second-order setting. In Section 3 we show that the same is true for a wide range of compactly supported kernels, including the biweight, provided our alternative definition of the critical bandwidth is used. In this sense the effect of nonmonotonicity of number of modes occurs at a relatively low level, and is not so great as to hinder the main features of a kernel density estimator. The mathematical argument behind this result is nonstandard, since a conventional approach relies on monotonicity, but nevertheless the result can be viewed as an extension of its counterpart for a Gaussian kernel.

All our methods and results have application to problems involving nonparametric regression, where only minor modifications are necessary. We have chosen to state them in the context of density estimation since passing in the reverse direction, from regression to density estimation, is not so straightforward; see the discussion by Chaudhuri and Marron [(2000), page 213].

There is no problem extending our results to the case of modes in estimators of density derivatives. As far as bimodality, or multimodality, is concerned, the main issue of interest is whether the bandwidth above which monotonicity of the number of modes (as a function of bandwidth) occurs is one for which the density estimator is multimodal with an appropriate number of modes. Indeed, if $k \geq 2$ is given then it is possible, when using a compactly supported kernel, that the density estimator will not have at least $k$ modes for a bandwidth, $h$, in a range $[h_0, \infty)$ where the number of modes is monotone in $h$. (This is relatively likely to occur if the actual density has strictly fewer than $k$ modes.) This possibility does not arise when $k = 1$, and for general $k$ it does not occur when using a Gaussian kernel. As a result, it is relatively unattractive to use compactly supported kernels in problems where strict multimodality is being investigated.

There is a diverse and extensive literature on bump hunting in nonparametric density estimation, much of it starting from contributions of Good and Gaskins (1980) and Silverman (1981). Formal and informal approaches to assessing modality include those of Hartigan and Hartigan (1985), Izenman and Sommer (1988), Roeder (1990, 1994), Cuevas and Gonzáles-Manteiga (1991), Müller and Sawitzki (1991), Minnotte and Scott (1993), Fisher, Mammen and Marron (1994), Escobar and West (1995), Polonik (1995a, b), Minnotte (1997), Chaudhuri and Marron (1999, 2000), Cheng and Hall (1999) and Fisher and Marron (2001). A small number of techniques, for example, the recent scale-space methods introduced by Chaudhuri and Marron (1999, 2000), rely on monotonicity of number of modes (as a function of bandwidth) in order to convey information. However, others, in particular formal or informal hypothesis testing approaches, require little more than the notion of a critical bandwidth and therefore suffer hardly at all from nonmonotonicity; as we show, lack of monotonicity occurs only for relatively small bandwidths. In these cases, and others, our results indicate that nonmonotonicity for popular kernels such as the biweight is generally not a significant problem. This serves to encourage their use in bump hunting problems.

## 2. Theory describing nonmonotonicity for non-Gaussian kernels.

2.1. *Preliminaries.* We say that a continuous density $f$ (or density estimator $\hat{f}$), continuously differentiable on its support, has just $k$ modes if it has only a finite number of points of inflection on its support, and just $k$ local maxima $x_1, \ldots, x_k$. The values of $x_j$ are called the modes of $f$. We say that $f$ is strictly unimodal if $f$ has just one mode in the sense defined above.

The assumption of continuity is made solely to simplify the definition of a density with $k$ modes; it may be weakened. Likewise we may remove the condition that the density has only isolated points of inflection on its support, although it should be appreciated that this alters the type of information contained in our results. We are not aware of any kernel used in practice that violates this condition.

Given a kernel $K$, bandwidth $h$ and sample $\mathcal{X} = \{X_1, \ldots, X_n\}$, let $\hat{f} = \hat{f}_h$, defined by

$$\hat{f}_h(x) = \frac{1}{nh} \sum_{i=1}^{n} K\left(\frac{x - X_i}{h}\right), \tag{2.1}$$

denote a conventional kernel density estimator. It is clear that if $K$ is strictly unimodal, continuous on the real line, and supported on a compact interval, and if the data $X_i$ are distinct, then for all sufficiently small $h$, $\hat{f}_h$ has just $n$ modes. We shall say that "the number of modes of the kernel estimator $\hat{f}_h$ is not monotone in $h$" if there exist $0 < h_1 < h_2$ such that the number of modes of $\hat{f}_{h_1}$ is strictly less than the number of modes of $\hat{f}_{h_2}$.

2.2. *Monotonicity of number of modes for large bandwidths.* First we introduce a unimodality condition:

(2.2)    $K$ is compactly supported and strictly unimodal, and is concave in a neighborhood of its mode.

Theorem 2.1 shows that (2.2) ensures $\hat{f}_h$ is unimodal for all sufficiently large $h$.

THEOREM 2.1. *If* (2.2) *holds, and if the data* $\mathcal{X}$ *come from a continuous distribution, then with probability* 1 *there exists a bandwidth* $\hat{h} = \hat{h}(\mathcal{X})$ *such that* $\hat{f}_h$ *is strictly unimodal for all* $h > \hat{h}$.

Proofs of Theorems 2.1 and 2.2 are given in Sections 5.1 and 5.2, respectively. A derivation of Theorem 2.3 is similar. Theorem 2.1 implies that the "critical bandwidth," given by

(2.3) $$h_{\text{crit}} = \inf\{h_1 > 0 \colon \hat{f}_h \text{ is unimodal for all } h > h_1\},$$

is well defined with probability 1. Moreover, assuming the sampled distribution is continuous, $P(h_{\text{crit}} > 0) = 1$. Throughout the paper, $h_{\text{crit}}$ is given by (2.3). Minnotte and Scott's (1993) mode tree algorithm permits calculation of $h_{\text{crit}}$. Without the algorithm, checking large bandwidths to see if the corresponding density estimator was unimodal could be computationally difficult.

2.3. *Theorems applicable to multiweight kernels.* Consider the condition

(2.4)    $K$ is a symmetric and strictly unimodal probability density with support equal to $\mathcal{I} = [-1, 1]$, continuous on the real line and continuously differentiable on $\mathcal{I}$, has two continuous derivatives in $[1 - \varepsilon, 1]$ for some $\varepsilon > 0$, and satisfies $K''(x) < 0$ for some $x \in (\frac{1}{2}, 1)$.

THEOREM 2.2. *If* (2.4) *holds, if* $n \geq 2$, *and if* $\mathcal{X} = \{X_1, \ldots, X_n\}$ *denotes a random sample drawn from a continuous distribution, then with probability* 1 *the number of modes of the kernel estimator* $\hat{f}_h$ *is not monotone in* $h$.

Any kernel of the form $K(x) = C_\theta (1 - x^2)^\theta$ on $\mathcal{I}$, where $0 < \theta < 5/2$ and $C_\theta$ ensures $\int K = 1$, satisfies (2.4). This class includes the uniweight (or Epanechnikov) and biweight kernels.

Theorem 2.2 does not address the triweight case ($\theta = 3$). In fact, when $n = 2$, and when $K$ is the triweight kernel and the sampled distribution is continuous, with probability 1 the number of modes of $\hat{f}_h$ is monotone nonincreasing as a function of the bandwidth. This result is available for more general kernels, too; it is sufficient that (2.2) hold and that $K$ be a symmetric probability density with support equal to $\mathcal{I}$, continuous on the real line, twice continuously differentiable on $\mathcal{I}$, with a unique point of inflection ($\xi$, say) on $(0, 1)$, and such that the only solutions $0 < x_1 \leq x_2 < 1$ of the equations $K'(x_1) = K'(x_2)$ and $K''(x_1) = -K''(x_2)$ are $x_1 = x_2 = \xi$. We ask too that $K$ have $2k \geq 4$ derivatives in a neighborhood of $\xi$, with $K^{(2j)}(\xi) = 0$ for $1 \leq j \leq k - 1$ and $K^{(2k)}(\xi) < 0$. These conditions hold with $k = 3$ when $K$ is the triweight kernel.

Our next result will show, however, that nonmonotonicity can occur with the triweight kernel provided $n \geq 3$. To this end, put $\kappa_\xi(x) = K(\xi + x) + K(\xi - x) + K(x)$, and assume that:

(2.5) $K$ is a symmetric and strictly unimodal probability density with support equal to $\mathcal{I} = [-1, 1]$, continuous on the real line, four times continuously differentiable on $\mathcal{I}$, and with the property that $\kappa''_\xi(0) > 0$, $\kappa'_\xi(\eta) = 0$ and $\kappa''_\xi(\eta) > 0$ for some $\xi, \eta \in (0, 1)$.

THEOREM 2.3. *If* (2.5) *holds, if* $n \geq 3$, *and if* $\mathcal{X} = \{X_1, \ldots, X_n\}$ *denotes a random sample drawn from a continuous distribution, then with strictly positive probability the number of modes of $\hat{f}_h$ is not monotone in $h$.*

Any kernel of the form $K(x) = C_\theta(1 - x^2)^\theta$ on $\mathcal{I}$, where $5/2 \leq \theta \leq 11/2$, satisfies (2.5). This class includes the triweight kernel, for which $\theta = 3$ and appropriate values of $\xi$ and $\eta$ are $\xi = 0.9$ and $\eta = 0.45$.

## 3. Critical bandwidths and bootstrap tests.

3.1. *Methodology.* The "classic" form of Silverman's (1981) bandwidth test for unimodality is based on computing a critical bandwidth that, in some sense, is as small as possible subject to the density estimator $\hat{f}_h$ at (2.1) being unimodal. If $K$ is a Gaussian density then there can be no ambiguity in defining the critical bandwidth: the number of modes is a monotone nonincreasing function of bandwidth, and so for any given dataset there is a bandwidth below which all density estimators have at least two modes, and above which all density estimators are unimodal [Schoenberg (1950)].

While this is not generally true for non-Gaussian kernels, that does not inhibit the definition of critical bandwidth given at (2.3). From a practical viewpoint it is quite feasible to define $h_{\text{crit}}$, as we do at (2.3), by decreasing through bandwidths for which $\hat{f}_h$ is unimodal, although it is generally not possible to define a critical bandwidth by increasing through bandwidths for which $\hat{f}_h$ is multimodal.

Silverman's (1981) bandwidth test for unimodality consists of rejecting the null hypothesis of unimodality if $h_{\mathrm{crit}}$ is "too large," where the latter is determined using the bootstrap. Specifically, put $\hat{f}_{\mathrm{crit}} = \hat{f}_{h_{\mathrm{crit}}}$, let $X_1^*, \ldots, X_n^*$ be a resample drawn by sampling randomly (conditional on $\mathcal{X}$) and with replacement from the distribution with density $\hat{f}_{\mathrm{crit}}$, and define

$$(3.1) \qquad \hat{f}_h^*(x) = (nh)^{-1} \sum_{i=1}^n K\left(\frac{x - X_i^*}{h}\right).$$

Let $h_{\mathrm{crit}}^*$ denote the version of $h_{\mathrm{crit}}$ in this setting, with $\hat{f}_h^*$ replacing $\hat{f}_h$ in the definition of $h_{\mathrm{crit}}$. Given a nominal level $\alpha$ for the test, the null hypothesis of unimodality is rejected if

$$(3.2) \qquad P(h_{\mathrm{crit}}^*/h_{\mathrm{crit}} \leq 1 | \mathcal{X}) \geq 1 - \alpha.$$

The technique, using our definition of $h_{\mathrm{crit}}$, can also be applied to assess unimodality in a subinterval of the support of a density.

3.2. *Large-sample properties of critical bandwidth.* Assume that:

(3.3) $f$ has two continuous derivatives when considered as a function restricted to its support, which we take to equal $\mathcal{S} = [a, b]$ where $-\infty < a < b < \infty$; that $f(a) = f(b) = 0$, $f'(a+) > 0$ and $f'(b-) < 0$; and that in the interior of $\mathcal{S}$ the equation $f'(x) = 0$ has a unique solution $x_0 \in (a, b)$, and that $f''(x_0) < 0$.

Assume too that:

(3.4) $K$ is a symmetric and strictly unimodal probability density with support equal to $\mathcal{I} = [-1, 1]$, is continuously differentiable on the real line, and has three bounded derivatives when viewed as a function defined only on $\mathcal{I}$.

This condition is satisfied by the biweight kernel, for example.

The part of condition (3.3) which asserts that $f$ decreases steeply to zero at either end of its support serves only to remove the effects of spurious "wiggles" in the tails of $\hat{f}_h$. Without such a constraint the size of the critical bandwidth can be determined by random clusters of data in the tails of the distribution. In practice such effects are usually excluded by restricting attention to the body of the distribution when formally testing for unimodality. However, there is a wide variety of ways of doing this, and for our purposes it is more appropriate to impose a condition which simply excludes tail effects. See Mammen, Marron and Fisher (1992) and Silverman (1983), for further discussion of this issue; they impose a condition close to (3.3).

Given a standard Brownian motion $W$, define the stochastic process

$$\omega(t, u) = f(x_0)^{1/2} u^{-1} \int K\left(\frac{t + v}{u}\right) dW(v) - \frac{1}{2}|f''(x_0)|t^2$$

for $-\infty < t < \infty$ and $u > 0$. Put $\omega'(t,u) = (\partial/\partial t)\omega(t,u)$. Our first result argues that for all sufficiently large $u$, the stochastic process $\omega(\cdot, u)$ is unimodal.

THEOREM 3.1. *Assume* (3.4) *holds. Then the probability that, for all $u \geq v$, $\omega(\cdot, u)$ has a unique local maximum and no local minimum [equivalently, $\omega'(\cdot, u)$ has a unique downcrossing of $0$ and no upcrossing of $0$] on the real line converges to $1$ as $v \to \infty$.*

It is readily shown that the probability that for some $u \geq v$, $\omega(\cdot, u)$ has both a local maximum and a local minimum, converges to $1$ as $v \downarrow 0$. This result and Theorem 3.1 imply that with probability 1 the infimum, $U_{\text{crit}}$ say, of the set of values $v > 0$ such that, for all $u \geq v$, $\omega(\cdot, u)$ has a unique local maximum and no local minimum, is well defined and strictly positive.

Our next result shows that $h_{\text{crit}}$ is asymptotically of conventional size $n^{-1/5}$, and that the "constant" of proportionality equals $U_{\text{crit}}$.

THEOREM 3.2. *Assume* (3.3) *and* (3.4) *hold. Then with probability 1, $h_{\text{crit}}$ is well defined for all sufficiently large $n$, and $n^{1/5}h_{\text{crit}}$ converges in distribution to $U_{\text{crit}}$.*

3.3. *Properties of bootstrap test.* First we describe large-sample properties of the distribution of $h^*_{\text{crit}}$, conditional on the data. Theorem 3.1 implies the existence of a unique point, $t = T$ say, at which

$$f(x_0)^{1/2} U_{\text{crit}}^{-2} \int K'\left(\frac{t+v}{U_{\text{crit}}}\right) dW(v) - |f''(x_0)|t$$

changes sign. Let $W^*$ denote a standard Wiener process independent of $W$, and put

$$\Omega^*(t,u) = f(x_0)^{1/2} u^{-2} \int K'\left(\frac{t+v}{u}\right) dW^*(v)$$
$$+ f(x_0)^{1/2} U_{\text{crit}}^{-2} \int K'\left(\frac{T+tu+v}{U_{\text{crit}}}\right) dW(v) - |f''(x_0)|(T+tu).$$

The argument used to prove Theorem 3.1 may be employed to show that the infimum $U^*_{\text{crit}}$ of the set of $v > 0$ such that, for all $u \geq v$, $\Omega^*(\cdot, u)$ has a unique downcrossing of $0$ and no upcrossing of $0$, is well defined and strictly positive. The strong approximation argument leading to the proof of Theorem 3.2 may be used to prove that, assuming both (3.3) and (3.4), and employing suitable constructions of $W$ and $W^*$,

$$\sup_{-\infty < x < \infty} |P(h^*_{\text{crit}}/h_{\text{crit}} \leq x | \mathcal{X}) - P(U^*_{\text{crit}}/U_{\text{crit}} \leq x | W)| \to 0$$

in probability.

It follows that the asymptotic level of the test defined at (3.2) is

(3.5) $$\pi(\alpha) = P\{P(U^*_{\text{crit}}/U_{\text{crit}} \leq 1|W) \geq 1-\alpha\}.$$

Note that $0 < \pi(\alpha) < 1$ for each $\alpha$. On the other hand, if the sampled density is not unimodal then $P(h^*_{\text{crit}}/h_{\text{crit}} \leq x) \to 1$ for each $x > 0$, and so the probability at (3.2) converges to 1. That is, when the null hypothesis is false, the probability that the test leads to rejection converges to 1 as $n \to \infty$. It is not true that $\pi(\alpha) = \alpha$, and this equality also fails in the case of a Gaussian kernel; see Hall and York (2001) for discussion of the size of the error.

## 4. Numerical properties for non-Gaussian kernels.

4.1. *Distribution of characteristics of nonmonotonicity.* The theoretical results in Section 2 may be illustrated using the mode tree of Minnotte and Scott (1993), for small samples and various kernels. In particular, the case $n = 3$ is treated in Figure 1. There we took the sample to be $\{X_1, X_2, X_3\} = \{-1, 0, 1\}$. We used $\sigma = 1/3$ for the Gaussian kernel. This gave effective support similar to that for the compact kernels.

The four panels in Figure 1 correspond, respectively, to the kernels: (a) Epanechnikov, (b) biweight, (c) triweight and (d) Gaussian. Panels (a)–(c) show that false modes appear at the points $\pm 1/2$ for each of the non-Gaussian kernels. This leads to nonmonotone behavior in each instance, and in fact as $h$ increases, 3 modes $\to 5 \to 3 \to 1$ for the Epanechnikov, $3 \to 5 \to 2 \to 3 \to 1$ for the biweight and $3 \to 4 \to 2 \to 1$ for the triweight. Clearly, the possibility of nonmonotonicity is very real for these kernels.

To further investigate the mode behavior of multiweight kernels for this three-point dataset, the number of modes was found for 500 values of $h$ and 480 choices of $\theta$ in the kernel $K_\theta(x) = C_\theta(1-x^2)^\theta$, ranging from 0.025 to 12. The result in "mode space" may be seen in Figure 2. The number of modes, between 1 and 6, is represented by the increasing density of six greyscale levels, as follows: 1 mode is indicated by the light grey in the north–west corner of the figure; 2 modes by the slightly darker adjacent region to its right, not touching any of the figure boundaries; 3 modes by the medium grey region that covers most of the south–east half of the figure, and also by the small area against the left-hand figure boundary immediately below the 1-mode region; 4 modes by the small sliver of a region between the 2-mode and 3-mode areas; 5 modes by the very dark patch which meets the left-hand figure boundary at values of $h$ between about 0.5 and 1.0; and a very small region of black, hardly detectable on the figure, representing 6 modes near $(\theta, h) = (2.5, 1.02)$.

The possibility of finding 6 modes in a density estimate from 3 points is demonstrated in Figure 3. Panel (a) shows a portion of the mode tree for the case $\theta = 2.5$, while panels (b) and (c) show the density estimate, in full and in modal close-up, respectively, for the estimate with $h = 1.02$ and the same kernel. The estimate is nearly flat, but six modes appear, ranging from small to extremely small.

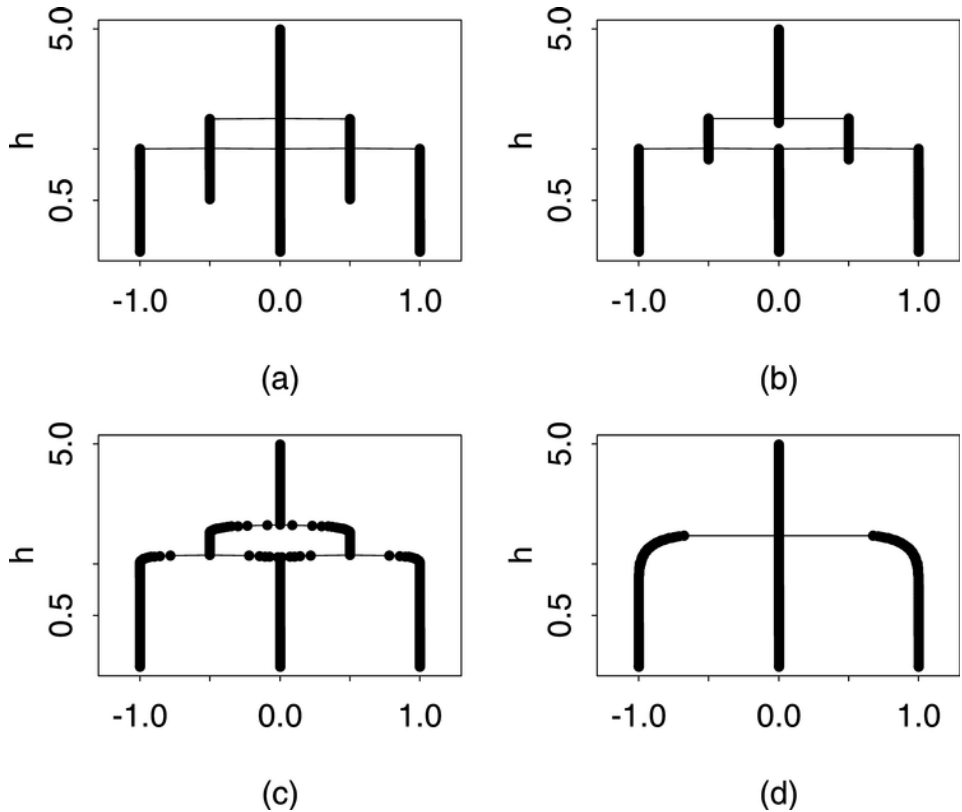

FIG. 1. *Mode trees for the $n = 3$ sample $\{-1, 0, 1\}$. Kernels used to produce the results in panels* (a)–(d) *are, respectively,* (a) *Epanechnikov,* (b) *biweight,* (c) *triweight and* (d) *Gaussian.*

Although clearly Figure 2 does not generalize directly to other datasets, it demonstrates both the complexity of the data-$\theta$-$h$ interactions with respect to modes, and the ubiquity of modal nonmonotonicity. Even though it is often assumed that $K_\theta$ provides a good approximation to the normal kernel for moderate values of $\theta$, the monotonicity property does not appear for this simple dataset until $\theta$ is close to 11.

Next we investigated the relationship between $h_{\text{crit}}$, defined in Section 2, and the bandwidth $h_{\text{nonm}}$, defined to be the smallest bandwidth at which nonmonotonicity appears as $h$ is decreased from $h_{\text{crit}}$. We drew 1000 samples of size $n$ from the distribution whose density was the Epanechnikov kernel, this choice being made because there the density estimates suffer in only minor ways from spurious modes in the tails. For each sample we computed density estimates using (a) Epanechnikov, (b) biweight and (c) triweight kernels, and formed the ratio $R = h_{\text{crit}}/h_{\text{nonm}}$. Estimates of the probability densities of $\log(R)$ are plotted in Figure 4, for $n = 10$ (dotted line), 100 (dashed line) and 1000 (solid line). For each estimate, the Gaussian kernel and the Sheather and Jones (1991) direct plug-in bandwidth were used. Note that both scales for the three panels vary considerably.

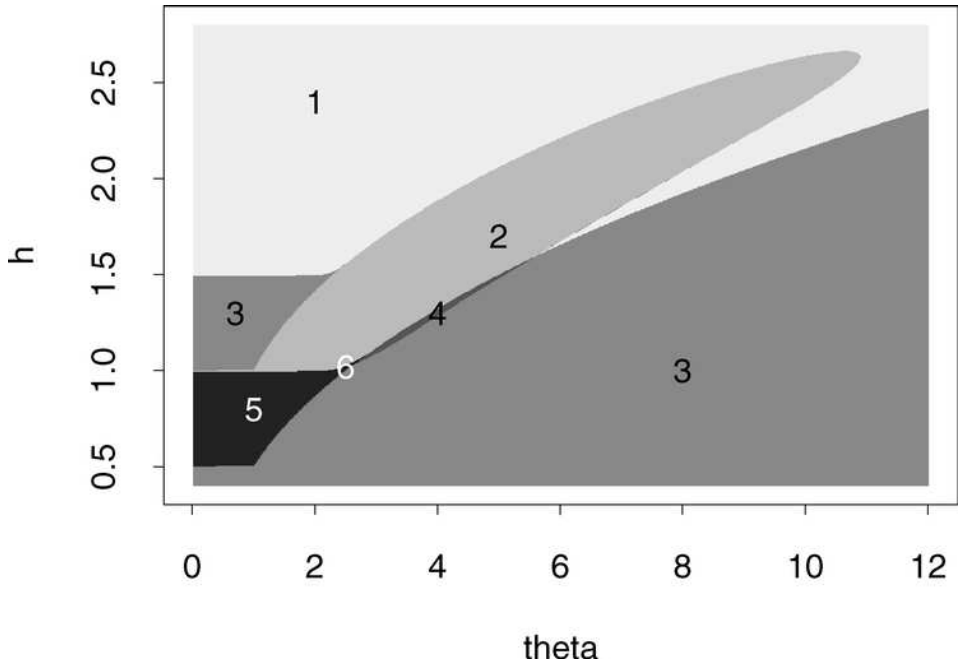

Fig. 2. *Image of "mode space." The figure shows the numbers of modes of density estimates computed from the sample $\{-1, 0, 1\}$, using the kernel $K_\theta(x) = C_\theta(1-x^2)^\theta$. Values of $\theta$ and $h$ are indicated on the horizontal and vertical axes, respectively. Mode counts range from 1 (light grey, upper left part of the figure), through 3 (medium grey, lower right part), to 6 (black).*

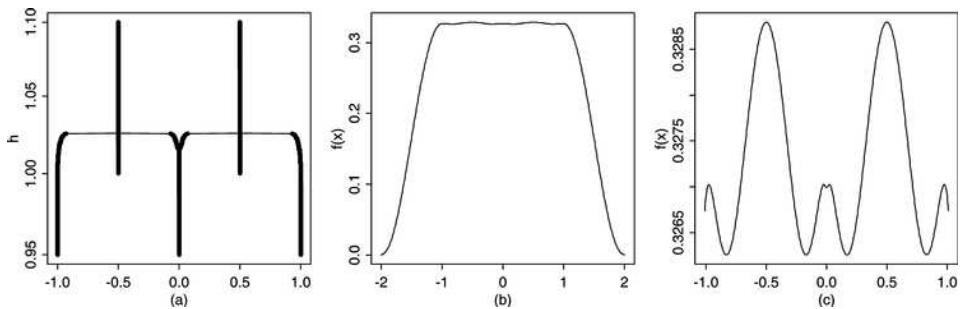

Fig. 3. *Six modes for estimates from the sample $\{-1, 0, 1\}$ using $K_{5/2}(x) = C_{5/2}(1-x^2)^{5/2}$. Panel (a) depicts part of the mode tree, while panel (b) shows the estimate using $h = 1.02$. Panel (c) displays a close-up of the modes from the same estimate.*

Panel (a) of Figure 4, for the case of the Epanechnikov kernel, shows that for all three sample sizes, nonmonotonicity tends to occur at a bandwidth that is close to $h_{\text{crit}}$. The biweight-kernel results presented in panel (b) show that nonmonotonicities are still very common, but that they now often appear at a bandwidth which is significantly smaller than $h_{\text{crit}}$. By way of contrast, panel (c) reveals that the triweight kernel is much less susceptible to nonmonotonicity, and that in this case sample size plays a larger role. The large peaks on the right in all three triweight estimates appear to be artifacts due to the discretized nature of the original estimates (400 points on $[-1, 1]$). It appears possible that a triweight kernel-based estimate suffers relatively few effective nonmonotonicities.

4.2. *Bump hunting.* In this section we summarize numerical information about the extent to which level accuracy of Silverman's bandwidth test, discussed in Section 3, is influenced by kernel type. Epanechnikov, biweight, triweight and Gaussian kernels are treated. It is well known that the Gaussian kernel produces asymptotically conservative tests, in the sense that the asymptotic level $\pi(\alpha)$, defined at (3.5), tends to be less than $\alpha$. It is of interest to learn what happens for other kernels.

Figure 5 illustrates results in the case of data simulated from the Beta$(3,4)$ distribution. The value of $h_{\text{crit}}$ was found by grid search. The bootstrap form, $h_{\text{crit}}^*$, of $h_{\text{crit}}$ was calculated by averaging over 500 bootstrap resamples from each sample, and the value of $\pi(\alpha)$ was approximated by averaging results over 100 replicates. The resulting curve approximations were slightly smoothed to reduce variability.

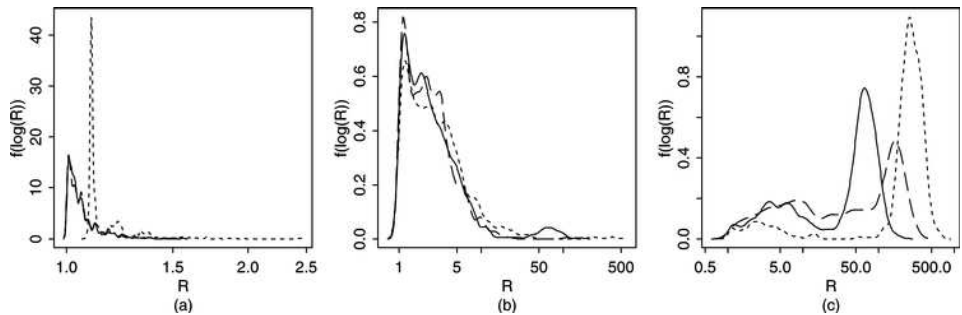

FIG. 4. *Estimates of probability density of* $\log(R = h_{\text{crit}}/h_{\text{nonm}})$. *Panels* (a)–(c) *correspond to the* (a) *Epanechnikov,* (b) *biweight and* (c) *triweight kernels, respectively. Sample sizes were* $n = 10$ *(represented by the dotted line),* $n = 100$ *(dashed line), and* $n = 1000$ *(solid line).*

Panels (a) and (b) in Figure 5 correspond to $n = 100$ and $n = 10000$, respectively. Each panel displays four approximations to $\pi(\alpha)$, indicated on the vertical axis, as functions of $\alpha$. The four curves represent the Epanechnikov (unbroken line), biweight (dotted line), triweight (dot-dash line) and Gaussian (dashed line) kernels. The conservative nature of the bandwidth test is indicated by the fact that each curve lies below the diagonal, with little to distinguish the different kernels. This lends support to the view that using non-Gaussian kernels when testing for modality does not substantially alter the conclusions of a test. The conservatism could be alleviated by using any of several available corrections, for example, that suggested by Hall and York (2001).

## 5. Technical arguments.

5.1. *Proof of Theorem* 2.1. Without loss of generality the mode of $K$ equals 0. It suffices to show that for each $n$ there exists $\varepsilon = \varepsilon(n) > 0$ such that, whenever $X_1, \ldots, X_n$ come from a continuous distribution supported on $[-\varepsilon, \varepsilon]$ [call this assumption (A)], the mixture density $g = n^{-1} \sum_i K(x + X_i)$ is strictly unimodal with probability 1.

If (A) holds then $g'(x) \geq 0$ whenever $x < -\varepsilon$, and equality occurs if and only if $K'(x + X_i) = 0$ for each $i$, which by assumption is true only at points of inflection of $K(\cdot + X_i)$ (by assumption there is only a finite number of these) or at points outside the support of $K(\cdot + X_i)$. Therefore if (A) holds then $g' \geq 0$ on $(-\infty, -\varepsilon)$ and $g' \leq 0$ on $(\varepsilon, \infty)$, with equality holding in either case only outside the support of $g$ or at points of inflection inside the support of $g$, there being at most a finite number of these. Call this property (P).

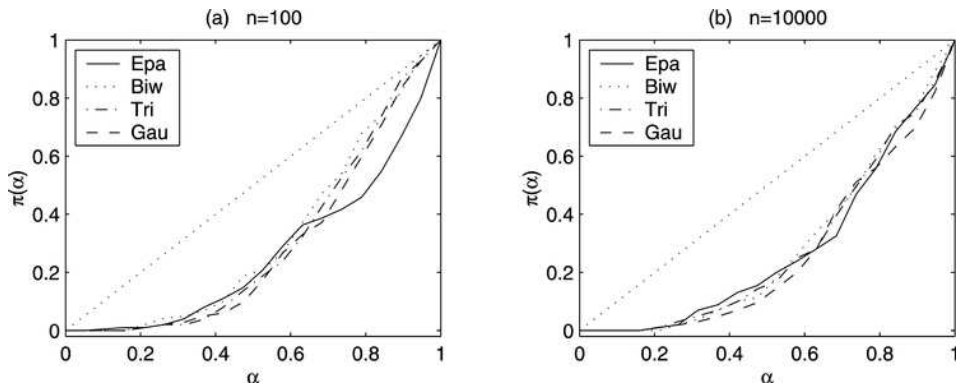

FIG. 5. *Level accuracy of bandwidth test. The four curves in each panel represent numerical approximations to levels of the bandwidth test when the Epanechnikov, biweight, triweight or Gaussian kernel is used to implement the test. Line types are as indicated in boxes. Panels* (a) *and* (b) *are for* $n = 100$ *and* $n = 10000$, *respectively.*

Since, for some $\eta > 0$, $K$ is concave in the neighborhood $(m - \eta, m + \eta)$ of $m$, then, provided assumption (A) holds for sufficiently small $\varepsilon$, $g$ is concave in $(m - \frac{1}{2}\eta, m + \frac{1}{2}\eta)$. Combining this property with (P) we deduce that $g$ is strictly unimodal on its support, except for the possibility that the set of points that gives a maximum of $g$ form a nondegenerate interval. However, this entails $\sum_i K'(x + X_i) = 0$ for all $x$ in that interval, which, since the sampled distribution is continuous and $K$ is strictly unimodal, holds with probability 0.

5.2. *Proof of Theorem* 2.2. Let $0 < \varepsilon < \frac{1}{2}$. Now,

$$g_\varepsilon(x) \equiv \tfrac{1}{2}\{K(-1 + \varepsilon + x) + K(1 - \varepsilon + x)\}$$
$$= K(1 - \varepsilon) + \tfrac{1}{2}x^2 K''(1 - \varepsilon) + o(x^2)$$

as $x \to 0$. Therefore, if $0 < \varepsilon < \frac{1}{2}$ then $K(-1 + \varepsilon + x) + K(1 - \varepsilon + x)$ is strictly concave in the neighborhood of the origin. It follows that the density $g_\varepsilon$ has at least three modes.

Equivalently, the density

$$(5.1) \qquad \frac{1}{2h} K\left(\frac{x - X_1}{h}\right) + \frac{1}{2h} K\left(\frac{x - X_2}{h}\right),$$

equal to the kernel density estimator computed from the sample $\{X_1, X_2\}$ of size $n = 2$, has at least three modes if $\frac{1}{2}|X_1 - X_2| < h < |X_1 - X_2|$, and has precisely two modes if $h \leq \frac{1}{2}|X_1 - X_2|$.

More generally, given a sample $\mathcal{X}$ of size $n \geq 2$ we may order the data as $X_{(1)} \leq \cdots \leq X_{(n)}$. Let $S_i = X_{(i+1)} - X_{(i)}$, for $1 \leq i \leq n - 1$, denote the $i$th spacing. If the sampled distribution is continuous then with probability 1 no two spacings are equal, and so they may be ranked in order of strictly increasing size, without ties. Let $S_{\min}$ denote the smallest spacing. Then with probability 1 the density estimator $\hat{f}_h$ has at least $n + 1$ modes if $\frac{1}{2}S_{\min} < h < S_{\min}$, and has precisely $n$ modes if $h \leq \frac{1}{2}S_{\min}$.

5.3. *Proof of Theorem* 3.1. Assumptions (3.4) imply that $\omega(t, u)$ has two continuous derivatives with respect to $t$, and that $u^{-1}\omega'(tu, u)$ is proportional to $\omega_1(t, u) = u^{-3/2}\omega_2(t, u) - ct$, where

$$\omega_2(t, u) = \int K''(t + v) W_u(v) \, dv,$$

$c = |f''(x_0)|/f(x_0)^{1/2}$ and $W_u(v) = -W(uv)/u^{1/2}$ is a standard Brownian motion. If $-\frac{1}{2} < t_1 < t_2 < \frac{1}{2}$ then

$$|\omega_2(t_1, u) - \omega_2(t_2, u)| \leq \int_{-1-t_1}^{1-t_2} |K''(t_2 + v) - K''(t_1 + v)||W_u(v)| \, dv$$

$$\begin{aligned}
&+ \int_{-1-t_2}^{-1-t_1} |K''(t_2+v)||W_u(v)|\,dv \\
&\qquad (5.2) \\
&+ \int_{1-t_2}^{1-t_1} |K''(t_1+v)||W_u(v)|\,dv \\
&\leq 4(t_2-t_1)\Big(\sup_{\mathcal{I}}|K''| + \sup_{\mathcal{I}}|K'''|\Big)S(u),
\end{aligned}$$

where $S(u) = \int_{-2\leq v\leq 2} |W_u(v)|\,dv$. (These bounds require $K$ to have three derivatives as a function on $\mathcal{I}$, but not as a function on the real line.) For each $\varepsilon > 0$, $S(u) = O(u^\varepsilon)$ with probability 1 as $u \to \infty$. Therefore, by (5.2),

$$(5.3) \qquad \sup_{-1 < t_1 < t_2 < 1} \left| \frac{\omega_2(t_1,u) - \omega_2(t_2,u)}{t_1 - t_2} \right| = O(u^\varepsilon)$$

with probability 1 as $u \to \infty$.

Solutions $t = \hat{t}$ of $\omega'(t,u) = 0$ are equivalently solutions of $\omega_1(t,u) = 0$, and may be shown by Taylor expansion to satisfy $\sup |\hat{t}| \to 0$ as $u \to \infty$, with probability 1, where the supremum is taken over all solutions. (It is straightforward to prove that with probability 1, at least one solution exists for all sufficiently large $u$.) Let $w > 0$ be given, and suppose the probability that for some $u \geq w$ at least two distinct solutions exist is bounded away from 0 (along a subsequence of values of $w$) as $w \to \infty$. Take $\hat{t}_1$ and $\hat{t}_2$ to be two such solutions, when $u \geq w$ and $w$ is an element of the subsequence. Then $u^{-3/2}\omega_2(\hat{t}_j,u) = c\hat{t}_j$ for each $j$, whence

$$(5.4) \qquad u^{-3/2}\frac{\omega_2(\hat{t}_1,u) - \omega_2(\hat{t}_2,u)}{\hat{t}_1 - \hat{t}_2} = c.$$

Result (5.3) implies, however, that with probability 1 the left-hand side of (5.4) converges to 0 as $u \to \infty$. On the other hand, the right-hand side is fixed and nonzero. This contradiction demonstrates the incorrectness of our assumption that two distinct solutions $\hat{t}_1$ and $\hat{t}_2$ of $\omega'(t,u) = 0$ exist, and proves the theorem.

5.4. *Proof of Theorem* 3.2. Put $h^0 = n^{-1/5}$, write $H = H(n)$ for a positive sequence such that $H(n) \to 0$ and $H(n)/h^0 \to \infty$, and redefine $h_{\text{crit}}$ to be the infimum of values $0 < h_1 \leq H(n)$ such that $\hat{f}_h$ is strongly unimodal for all $h \geq h_1$. It is readily proved that the probability that this version of $h_{\text{crit}}$, and the version defined at (2.3), are identical converges to 1 as $n \to \infty$. Therefore it is sufficient to prove that for the new version, $n^{1/5} h_{\text{crit}} \to U_{\text{crit}}$ in distribution.

Let $\widehat{\mathcal{S}}_h$ denote the support of $\hat{f}_h$, and write $\widehat{\mathcal{S}}_h^0$ for the interior of $\widehat{\mathcal{S}}_h$. Using strong approximation of the empirical distribution by a Brownian bridge [Komlós, Major and Tusnády (1976)], it may be shown that for each $C_1, \varepsilon > 0$ there exists $C_2 = C_2(C_1, \varepsilon) > 0$ such that for all sufficiently large $n$,

(5.5) $\quad P\{$for all $h \in [C_1 h^0, H(n)]$, $\hat{f}'_h(x) > 0$
for $x \in \widehat{\mathcal{S}}^0_h$ such that $x \leq x_0 - C_2 h$,
and $\hat{f}'_h(x) < 0$ for $x \in \widehat{\mathcal{S}}^0_h$ such that $x \geq x_0 + C_2 h\} \geq 1 - \varepsilon.$

The method of proof consists of showing first that for each $h$, $Ehf'_h(x)$ is strictly positive on $(a - h, x_0 - C_2 h)$ and strictly negative on $(x_0 + C_2 h, b + h)$, and thence demonstrating that for each $C_1, \varepsilon > 0$ there exists $C_2 = C_2(C_1, \varepsilon) > 0$ such that for all sufficiently large $n$,

$P\{$for all $h \in [C_1 h^0, H(n)]$ and all $x \in (a - h, b + h)$ for
which $|x - x_0| > C_2 h$, $|\hat{f}'_h(x) - Ehf'_h(x)| > \frac{1}{2}|Ehf'_h(x)|\} \geq 1 - \varepsilon.$

The same strong approximation methods may be used to prove that for each $C_2, \varepsilon > 0$ there exists $C_3 = C_3(C_2, \varepsilon) > 0$ such that for all sufficiently large $n$,

(5.6) $\quad P\left\{\sup_{h \in [C_3 h^0, H(n)]} \sup_{|x - x_0| \leq C_2 h} |\hat{f}''_h(x) - E\hat{f}''_h(x)| > \varepsilon\right\} \geq 1 - \varepsilon.$

[In each case the arguments are broadly similar to those of Mammen, Marron and Fisher (1992). See also Silverman (1983).]

Note too that for each $C_1 > 0$, $E\{\hat{f}''_h(x)\} = f''(x_0) + o(1)$ uniformly in $|x - x_0| \leq C_1 h$ and $h \leq H(n)$. Combining this result with (5.6) we deduce that for each $C_2, \varepsilon > 0$ there exists $C_3 = C_3(C_2, \varepsilon) > 0$ such that for all sufficiently large $n$,

(5.7) $\quad P\{$for all $h \in [C_3 h^0, H(n)]$,
$\hat{f}_h$ is strictly concave on $(x_0 - C_2 h, x_0 + C_2 h)\} \geq 1 - \varepsilon.$

Together (5.5) and (5.7) imply that for each $\varepsilon > 0$ there exists $C_3 = C_3(\varepsilon) > 0$ such that for all sufficiently large $n$,

(5.8) $\quad P\{$for all $h \in [C_3 h^0, H(n)]$,
$\hat{f}_h$ is strictly unimodal on its support$\} \geq 1 - \varepsilon.$

Strong approximation methods may also be used to prove the existence of a Brownian motion $W$ such that, defining $h_u = u h^0$,

$$A_j(t, u) = n^{(2-j)/5}\{\hat{f}^{(j)}_{h_u}(x_0 + h^0 t) - E\hat{f}^{(j)}_{h_u}(x_0 + h^0 t)\} \quad \text{and}$$

$$a_j(t, u) = f(x_0)^{1/2} u^{-(j+1)} \int K^{(j)}\left(\frac{v + t}{u}\right) dW(v),$$

where $j = 0$, 1 or 2, we have for each $0 < C_1 < C_3 < \infty$ and each $C_2 > 0$,

$$(5.9) \quad P\left\{\sup_{|t| \leq C_2} \sup_{C_1 \leq u \leq C_3} |A_j(t,u) - a_j(t,u)| \geq n^{-\delta}\right\} = O(n^{-\lambda})$$

for some $\delta > 0$ and all $C_2, \lambda > 0$. Observe too that if $j = 1, 2$,

$$(5.10) \quad n^{(2-j)/5}\{E\hat{f}_{h_u}^{(j)}(x_0 + h^0 t) - f^{(j)}(x_0)\} - tf''(x_0)I(j=1) \to 0$$

uniformly in $|t| \leq C_2$ and $C_1 \leq u \leq C_3$.

Denote by $N = N(C_2, u)$ the number of crossings of 0 made by the process $a_1(\cdot, u)$ in $[-C_2, C_2]$, and let the crossings be $T_1(u), \ldots, T_N$. For each $C_2 > 0$ and $0 < C_1 < C_3 < \infty$, the value of $\sup_{C_1 \leq u \leq C_3} N(C_2, u)$ is finite with probability 1. The continuous, nondegenerate property of the joint distributions of $a_1(\cdot, u)$ and $a_2(\cdot, u)$ implies that

$$\lim_{\varepsilon \to 0} P[|a_2\{T_i(u), u\}| > \varepsilon \text{ for } 1 \leq i \leq N(C_2, u) \text{ and } C_1 \leq u \leq C_3] = 1.$$

Note too that $\omega'(t, u) = a_1(t, u) + tf''(x_0)$. The results immediately above, and (5.9) and (5.10), imply that

$$(5.11) \quad \begin{aligned} &P\{\text{for each } C_1 \leq u \leq C_3, \text{ the number of downcrossings of 0 made} \\ &\text{by } \hat{f}'_{h_u}(x_0 + h^0 t) \text{ for } t \in [-C_2, C_2], \text{ equals the number} \\ &\text{of downcrossings of 0 made by } \omega'(t, u) \text{ for } t \in [-C_2, C_2]\} \to 1 \end{aligned}$$

as $n \to \infty$.

Analogously to (5.8), but more simply, it may be shown that for each $\varepsilon > 0$ there exists $C_3 = C_3(\varepsilon) > 0$ such that

$$(5.12) \quad \begin{aligned} &P\{\text{for all } u > C_3, \text{ there is a unique } t = \hat{t} \in (-\infty, \infty) \text{ at} \\ &\text{which } \omega'(t, u) \text{ vanishes, and } \hat{t} \text{ is a downcrossing}\} \geq 1 - \varepsilon. \end{aligned}$$

Combining (5.5), (5.8), (5.11) and (5.12), we deduce that for all $C_1 > 0$,

$P\{\text{for all } h \in [C_1 h^0, H(n)], \text{ the number of downcrossings of 0 made}$

$\text{by } \hat{f}_h \text{ equals the number of downcrossings of 0 made by } \omega'(t, u)\} \geq 1 - \varepsilon.$

The theorem follows from this result and (5.12).

**Acknowledgment.** We are grateful for the helpful and constructive comments of three reviewers.

## REFERENCES

Chaudhuri, P. and Marron, J. S. (1999). SiZer for exploration of structures in curves. *J. Amer. Statist. Assoc.* **94** 807–823. MR1723347

Chaudhuri, P. and Marron, J. S. (2000). Scale space view of curve estimation. *Ann. Statist.* **28** 408–428. MR1790003

Cheng, M.-Y. and Hall, P. (1999). Mode testing in difficult cases. *Ann. Statist.* **27** 1294–1315. MR1740110

Cuevas, A. and González-Manteiga, W. (1991). Data-driven smoothing based on convexity properties. In *Nonparametric Functional Estimation and Related Topics* (G. Roussas, ed.) 225–240. Kluwer, Dordrecht. MR1154331

Escobar, M. D. and West, M. (1995). Bayesian density estimation and inference using mixtures. *J. Amer. Statist. Assoc.* **90** 577–588. MR1340510

Fisher, N. I., Mammen, E. and Marron, J. S. (1994). Testing for multimodality. *Comput. Statist. Data Anal.* **18** 499–512. MR1310472

Fisher, N. I. and Marron, J. S. (2001). Mode testing via the excess mass estimate. *Biometrika* **88** 499–517. MR1844848

Good, I. J. and Gaskins, R. A. (1980). Density estimation and bump-hunting by the penalized likelihood method exemplified by scattering and meteorite data (with discussion). *J. Amer. Statist. Assoc.* **75** 42–73. MR568579

Hall, P. and York, M. (2001). On the calibration of Silverman's test for multimodality. *Statist. Sinica* **11** 515–536. MR1844538

Hartigan, J. A. and Hartigan, P. M. (1985). The DIP test of unimodality. *Ann. Statist.* **13** 70–84. MR773153

Izenman, A. J. and Sommer, C. (1988). Philatelic mixtures and multimodal densities. *J. Amer. Statist. Assoc.* **83** 941–953.

Komlós, J., Major, P. and Tusnády, G. (1976). An approximation of partial sums of independent rv's, and the sample df. II. *Z. Wahrsch. Verv. Gebiete* **34** 33–58. MR402883

Mammen, E., Marron, J. S. and Fisher, N. I. (1992). Some asymptotics for multimodality tests based on kernel density estimates. *Probab. Theory Related Fields* **91** 115–132. MR1142765

Minnotte, M. C. (1997). Nonparametric testing of the existence of modes. *Ann. Statist.* **25** 1646–1660. MR1463568

Minnotte, M. C. and Scott, D. W. (1993). The mode tree: A tool for visualization of nonparametric density estimates. *J. Comput. Graph. Statist.* **2** 51–68.

Müller, D. W and Sawitzki, G. (1991). Excess mass estimates and tests for multimodality. *J. Amer. Statist. Assoc.* **86** 738–746. MR1147099

Polonik, W. (1995a). Measuring mass concentrations and estimating density contour clusters—an excess mass approach. *Ann. Statist.* **23** 855–881. MR1345204

Polonik, W. (1995b). Density estimation under qualitative assumptions in higher dimensions. *J. Multivariate Anal.* **55** 61–81. MR1365636

Roeder, K. (1990). Density estimation with confidence sets exemplified by superclusters and voids in the galaxies. *J. Amer. Statist. Assoc.* **85** 617–624.

Roeder, K. (1994). A graphical technique for determining the number of components in a mixture of normals. *J. Amer. Statist. Assoc.* **89** 487–495. MR1294074

Schoenberg, I. J. (1950). On Pólya frequency functions. II. Variation-diminishing integral operators of the convolution type. *Acta Sci. Math. (Szeged)* **12** 97–106. MR35861

Sheather, S. J. and Jones, M. C. (1991). A reliable data-based bandwidth selection method for kernel density estimation. *J. Roy. Statist. Soc. Ser. B* **53** 683–690. MR1125725

Silverman, B. W. (1981). Using kernel density estimates to investigate multimodality. *J. Roy. Statist. Soc. Ser. B* **43** 97–99. MR610384

Silverman, B. W. (1983). Some properties of a test for multimodality based on kernel density estimates. In *Probability, Statistics and Analysis* (J. F. C. Kingman and G. E. H. Reuter, eds.) 248–259. Cambridge Univ. Press. MR696032

P. Hall  
Centre for Mathematics  
  and its Applications  
Australian National University  
Canberra, ACT 0200  
Australia  
e-mail: Peter.Hall@maths.anu.edu.au

M. C. Minnotte  
Department of Mathematics  
  and Statistics  
Utah State University  
Logan, Utah 84322-3900  
USA  
e-mail: minnotte@math.usu.edu

C. Zhang  
Department of Statistics  
University of Wisconsin  
1300 University Avenue  
Madison, Wisconsin 53706  
USA  
e-mail: cmzhang@stat.wisc.edu